\documentclass[10pt,a4paper]{amsart}
\usepackage{amssymb}
\usepackage{amsmath}
\usepackage{amsfonts}

\newcommand{\Q}{\mathbb{Q}}

\newcommand{\Z}{\mathbb{Z}}
\newcommand{\F}{\mathbb{F}}

\newtheorem{main-dummy}{Main-Dummy}
\newtheorem{dummy}{Dummy}

\numberwithin{dummy}{section}
\numberwithin{equation}{section}

\newtheorem{main-theorem}[main-dummy]{Theorem}
\newtheorem*{Lucas}{Lucas' Theorem}

\newtheorem{lemma}[dummy]{Lemma}
\newtheorem{theorem}[dummy]{Theorem}
\newtheorem{prop}[dummy]{Proposition}
\newtheorem{cor}[dummy]{Corollary}

\theoremstyle{definition}
\newtheorem{definition}[dummy]{Definition}

\theoremstyle{remark}
\newtheorem{rem}[dummy]{Remark}

\begin{document}
\bibliographystyle{amsalpha}
\author{Sandro Mattarei}

\email{mattarei@science.unitn.it}

\urladdr{http://www-math.science.unitn.it/\~{ }mattarei/}

\address{Dipartimento di Matematica\\
  Universit\`a degli Studi di Trento\\
  via Sommarive 14\\
  I-38050 Povo (Trento)\\
  Italy}

\title[Modular periodicity of binomial coefficients]{Modular periodicity of binomial coefficients}

\begin{abstract}
We prove that if the signed binomial coefficient
$(-1)^i\binom{k}{i}$ viewed modulo $p$ is a periodic function of
$i$ with period $h$ in the range $0\le i\le k$, then $k+1$ is a
power of $p$, provided $h$ is prime to $p$ and not too large compared to $k$.
(In particular, $2h\le k$ suffices.)
As an application, we prove that if $G$ and $H$ are multiplicative subgroups of a finite field,
with $H<G$, and such that
$1-\alpha\in G$ for all $\alpha\in G\setminus H$,
then $G\cup\{0\}$ is a subfield.
\end{abstract}




\thanks{The  author  is grateful  to  Ministero dell'Istruzione, dell'Universit\`a  e
  della  Ricerca, Italy,  for  financial  support of the
  project ``Graded Lie algebras  and pro-$p$-groups of finite width''.}

\maketitle

\thispagestyle{empty}

\section{Introduction}\label{sec:intro}

Binomial coefficients $\binom{k}{i}$
display several periodicity properties when viewed modulo a prime number $p$,
which are essentially related to the fact that $(a+b)^p=a^p+b^p$
in a commutative ring of characteristic $p$.
For example, for any (not necessarily positive) integer $i$ the function $k\mapsto\binom{k}{i}\pmod{p}$
is periodic (for $k\in\Z$) with period the smallest
power of $p$ which is greater than $i$.
(In fact, this is the shortest period, according to~\cite[Theorem~4.8]{Fray}, which, more generally,
gives the shortest period of the function
$k\mapsto\binom{k}{i}\pmod{p^r}$.)

If we take the identity
\[
(1+x)^k
=\sum_{i\in\Z}
\binom{k}{i}x^i
=\sum_{i\ge 0}
\binom{k}{i}x^i
\]
in the ring of formal power series $\Z[[x]]$ as the definition
of the integers $\binom{k}{i}$ for $k,i\in\Z$,
many properties of binomial coefficients modulo a prime $p$
(or {\em modular} properties for short)
can be proved conveniently by means of identities in the quotient ring
$\Z[[x]]/p\,\Z[[x]]\cong\F_p[[x]]$, where $\F_p$ is the field of $p$ elements
(or, more simply, the polynomial ring $\F_p[x]$ if we restrict to $k\ge 0$).
In particular, the periodicity property stated in the previous paragraph is an immediate consequence of the identity
\[
(1+x)^{k+p^r}=(1+x)^k(1+x^{p^r})=(1+x)^k+x^{p^r}(1+x)^k
\]
in $\F_p[[x]]$, where $p^r$ is the smallest power of $p$ exceeding $i$.
A natural generalization of this argument leads to what is known as Lucas' theorem~\cite{Lucas}.
For a given prime power $p^r$, we use the following notation:
an integer $b$ can be written uniquely as $b=b'p^r+b''$ where $0\le b''<p^r$.

\begin{Lucas}
$\binom{k}{i}\equiv\binom{k'}{i'}\binom{k''}{i''}\pmod{p}$.
\end{Lucas}

\begin{proof}
Apply the identity
$(1+x)^{k}=(1+x^{p^r})^{k'}(1+x)^{k''}$
in $\F_p[[x]]$.
\end{proof}

Recursive application of Lucas' theorem reduces the evaluation  modulo $p$ of arbitrary binomial coefficients
to that of binomial coefficients where both entries are less than $p$.

In this paper we explore the occurrences of periodicity with period $h$ prime to $p$.
Of course, the existence of such a periodicity in $\binom{k}{i}\pmod{p}$ with respect to the upper entry $k$
together with some periodicity of period $p^r$ (which we have just seen to exist as soon as
$p^r$ exceeds $i$) would imply that $1=(h,p^r)$ is also a period,
and one would easily conclude that $i\le 0$ (and hence $\binom{k}{i}=0$ for all $k$).
We allow, however, periodicity to occur only {\em in a certain range},
as in the following definition.

\begin{definition}\label{def:periodic}
A function $i\mapsto f(i)$ defined on a subset of the integers and taking values in any set
is periodic in the range $a\le i\le b$ with period $h$ (a positive integer)
if $f$ is defined in that range and $f(i+h)=f(i)$ whenever $a\le i\le b-h$.
\end{definition}

In this paper we consider periodicity of binomial coefficients with respect
to the lower entry $i$ rather than the upper entry $k$.
Since $\binom{k}{i}$ for $k\ge 0$ vanishes unless $0\le i\le k$,
periodicity for all $i\in\Z$ never occurs,
and the natural range to consider for Definition~\ref{def:periodic} is $0\le i\le k$.
We will see in Section~\ref{sec:variations} that modular periodicity in $k$ in
appropriate ranges is to some extent equivalent to that in $i$.
Furthermore, it will be more natural to consider {\em signed} binomial coefficients
$(-1)^i\binom{k}{i}$, which for a given $k$ are the coefficients of the formal power series
obtained by expanding $(1-x)^k$.
Some reasons for using these or other alternating signs will also be mentioned in Section~\ref{sec:variations}.

The paradigm of periodicity of $(-1)^i\binom{k}{i}\pmod{p}$ with respect to $i$ in the range
$0\le i\le k$ occurs when $k$ is one less than a power of $p$, in which case we have
$(-1)^i\binom{k}{i}\equiv 1\pmod{p}$ in the whole range,
as follows from
$\sum_{i=0}^k(-1)^i\binom{k}{i}=(1-x)^k=(1-x^{k+1})/(1-x)=\sum_{i=0}^k x^i$.
Our main result
asserts that this is the only occurrence of periodicity of
$(-1)^i\binom{k}{i}\pmod{p}$ with respect to $i$ in the range
$0\le i\le k$ with period $h$ prime to $p$, provided $h$ is not too large with respect to $k$.
Taking $h\le k/2$ is sufficient, but we give weaker and more precise assumptions in Theorem~\ref{thm:binomial} below.
We will show in Remark~\ref{rem:best-possible} that these assumptions cannot be weakened any further.
Theorem~\ref{thm:binomial} will be generalised in Corollary~\ref{cor:binomial}
by removing the hypothesis that $h$ is not a multiple of $p$
and suitably adapting the conclusion.

\begin{main-theorem}\label{thm:binomial}
Let $p$ be a prime, and let $h,k$ be positive integers with $p\nmid h$ and $k\ge 5$.
Suppose that $3h<2k+5$ or $p=3$ and $5h<4k+9$.
Suppose that the signed binomial coefficient $(-1)^i\binom{k}{i}$ viewed modulo $p$
is a periodic function of $i$ with period $h$ in the range $0\le i\le k$, that is,
\[
\binom{k}{i+h}\equiv
(-1)^h\binom{k}{i}\pmod{p}
\]
for $0\le i\le k-h$.
Then $k+1$ is a power of $p$.
\end{main-theorem}

The special case of Theorem~\ref{thm:binomial} where $h=1$ follows quickly from Lucas' Theorem,
but we give here a different proof which illustrates techniques that will be used later.
If $(-1)^i\binom{k}{i}\equiv 1\pmod{p}$ for all $0\le i\le k$, then
$(1-x)^k-1=\sum_{i=1}^{k}x^i=x(x^k-1)/(x-1)$ in $\F_p[x]$.
Note that $k$ must be prime to $p$ because $k=\binom{k}{1}\equiv -1\pmod{p}$.
Hence the group $G$ of $k$-th roots of unity in a splitting field for $x^k-1$ over $\F_p$ has order $k$.
The polynomial identity which we have found tells us that
$1-\alpha\in G\cup\{0\}$ for all $\alpha\in G\cup\{0\}$.
Since $G\cup\{0\}$ is multiplicatively closed, we infer that
$\beta-\alpha\beta=\beta(1-\alpha)\in G\cup\{0\}$ for all $\alpha,\beta\in G\cup\{0\}$.
Thus $G\cup\{0\}$ is also closed with respect to taking differences of its elements.
(For the special case of taking opposites write $-\gamma$ as $(\beta-\gamma)-\beta$.)
Hence $G\cup\{0\}$ is a field of order $k+1$, and we conclude that $k+1$ is a power of $p$.

Despite the elegance of the approach with polynomials and finite fields
in this special case, in more complex arguments
Lucas' theorem proves invaluable for evaluating binomial coefficients modulo $p$,
and is the basic tool on which our proof of Theorem~\ref{thm:binomial} rests.
Nevertheless, the argument which we have used above to deal with the special case $h=1$
will be reversed in Section~\ref{sec:finite-fields} to deduce from Theorem~\ref{thm:binomial} the following consequence
about multiplicative subgroups of finite fields.
\begin{main-theorem}\label{thm:near-field}
Let $G$ be a multiplicative subgroup of a finite field $\F_q$
which generates $\F_q$ as a field,
and let $H$ be a proper subgroup of $G$.
Suppose that
$1-\alpha\in G$ for all $\alpha\in G\setminus H$.
Then $G=\F_{q}^\ast$.
\end{main-theorem}

This result can also be proved directly using an argument of Leep and Shapiro
from~\cite{LeeSha}, as we explain in Section~\ref{sec:finite-fields}.
We also give an alternative approach to this problem
based on known bounds for the number of points of Fermat curves over finite fields.
In particular, this gives yet another proof of Theorem~\ref{thm:near-field}
in case $q$ is odd.

In Section~\ref{sec:variations} we present an analogue of Theorem~\ref{thm:binomial} for unsigned binomial
coefficients, and the corresponding analogue of Theorem~\ref{thm:near-field}.
Furthermore, we translate Theorem~\ref{thm:binomial} into a result concerning modular periodicity
of binomial coefficients with respect to the upper entry.

\medskip

I am grateful to the referee for several comments which led to an improved exposition.

\section{Modular periodicity of binomial coefficients}\label{sec:periodicity}

Our proof of Theorem~\ref{thm:binomial} relies on the following result.

\begin{prop}\label{prop:binomial}
Suppose $k$ is a positive integer and $k+1=p^rt$, where $p$ is prime and $p\nmid t$.
If $t>1$ then the signed binomial coefficient $(-1)^i\binom{k}{i}$ modulo $p$
cannot be a periodic function of $i$ in the range $0\le i\le k$
for any period $h\le k-p^r$ with $p\nmid h$.
\end{prop}

\begin{proof}
All congruences in this proof will be taken modulo $p$.
Suppose there is some period $h\le k-p^r$ with $p\nmid h$.
We start by proving that $k\equiv -1$.
Because of the periodicity hypothesis we have
$\binom{k}{h}\equiv (-1)^h$ and
$\binom{k}{h+1}\equiv (-1)^h\,k$.
The identity
$\binom{k}{h}(k-h)=\binom{k}{h+1}(h+1)$,
which follows from the factorial formula for binomial coefficients,
implies that
$k-h\equiv k(h+1)$.
Hence $(k+1)h\equiv 0$,
and we conclude that $k\equiv -1$ since $p\nmid h$.
This gives a contradiction if $r=0$, and hence we assume that $r>0$ from now on.

With notation as in Lucas' Theorem,
$k=k'p^r+k''$ where $k'=t-1$ and $k''=p^r-1$.
Periodicity at $i=0$ and Lucas' Theorem imply that
\[
(-1)^h=
(-1)^h
\binom{k}{0}\equiv
\binom{k}{h}\equiv
\binom{k'}{h'}\binom{p^r-1}{h''}\equiv
\binom{k'}{h'}(-1)^{h''}.
\]
Consider now periodicity at $i=p^r-1$, which holds because $p^r-1\le k-h$.
According to Lucas' Theorem we have
$\binom{k}{p^r-1}\equiv
\binom{k'}{0}\binom{p^r-1}{p^r-1}=1$.
Similarly, since $h+p^r-1=p^r(h'+1)+(h''-1)$,
with $h''>0$ because $p\nmid h$,
Lucas' Theorem implies that
$\binom{k}{h+p^r-1}\equiv
\binom{k'}{h'+1}\binom{p^r-1}{h''-1}=
\binom{k'}{h'+1}(-1)^{h''-1}$.
Thus, periodicity at $i=p^r-1$ gives us the equation
\[
(-1)^{h}=
(-1)^h\binom{k}{p^r-1}\equiv
\binom{k}{h+p^r-1}\equiv
\binom{k'}{h'+1}(-1)^{h''-1}.
\]
Combining this equation with that found earlier we obtain that
$\binom{k'}{h'+1}\equiv -\binom{k'}{h'}\not\equiv 0$.
The identity
$\binom{k'}{h'}(k'-h')=\binom{k'}{h'+1}(h'+1)$ implies that
$k'\equiv -1$.
But then $t=k'+1\equiv 0$, providing the desired contradiction.
\end{proof}

\begin{proof}[Proof of Theorem~\ref{thm:binomial}]
Suppose for a contradiction that $k+1$ is not a power of $p$, and let
$p^r$ be the highest power of $p$ which divides $k+1$.
Hence $k\ge 2p^r-1$ if $p$ is odd, and $k\ge 3p^r-1$ if $p=2$.
If $h+p^r\le k$ Proposition~\ref{prop:binomial} applies and yields a contradiction,
hence we may assume that $h+p^r\ge k+1$.
Then our hypotheses $3h<2k+5$, or $p=3$ and $5h<4k+9$, imply that
$k<3p^r+2$, or $p=3$ and $k<5p^r+4$.
It remains to check individually the cases
$k=2p^r-1$ for $p$ odd,
$k=3p^r-1$ for $p\not=3$, and
$k=4p^r-1$ or $5p^r-1$ for $p=3$.
(In addition to these cases, when $r\le 1$ the above inequalities also allow $k=5$ for $p=2$ and $k=6,7$ for $p=3$,
which can all be excluded by inspection.)

If $k=2p^r-1$ and $p$ is odd we have $(-1)^i\binom{k}{i}\equiv 1\pmod{p}$ for $0\le i<p^r$ and
$(-1)^i\binom{k}{i}\equiv -1\pmod{p}$ for $p^r\le i<2p^r$,
hence no periodicity occurs.
Similarly, if $k=4\cdot 3^r-1$ we have $(-1)^i\binom{k}{i}\equiv 1\pmod{3}$ for $0\le i<3^r$,
$(-1)^i\binom{k}{i}\equiv 0\pmod{3}$ for $3^r\le i<3\cdot 3^r$, and
$(-1)^i\binom{k}{i}\equiv -1\pmod{3}$ for $3\cdot 3^r\le i<4\cdot 3^r$;
hence no periodicity occurs.

If $k=3p^r-1$ and $p\not=3$ our hypothesis that $3h<2k+5$ implies that $h<2p^r+1$, hence $h\le 2p^r-1$ because $p\nmid h$,
and so $h+p^r\le k$ actually does hold, and Proposition~\ref{prop:binomial} applies.
Similarly, if $k=5\cdot 3^r-1$ and $p=3$ our hypothesis that $5h<4k+9$ implies that $h<4\cdot 3^r+1$,
hence $h\le 4\cdot 3^r-1$ because $3\nmid h$,
and so $h+3^r\le k$ holds and Proposition~\ref{prop:binomial} applies.
\end{proof}

The following statement is weaker than Theorem~\ref{thm:binomial} but easier to remember,
and is sufficient for some applications, such as Theorem~\ref{thm:near-field}, which we prove in the next section.

\begin{cor}\label{cor:binomial-weaker}
Let $p$ be a prime, and let $h,k$ be positive integers with $p\nmid h$ and $2h\le k$.
Suppose that the signed binomial coefficient $(-1)^i\binom{k}{i}$ viewed modulo $p$
is a periodic function of $i$ with period $h$ in the range $0\le i\le k$.
Then $k+1$ is a power of $p$.
\end{cor}

\begin{proof}
The conclusion follows from Theorem~\ref{thm:binomial} if $k\ge 5$, and by inspection in the remaining cases.
\end{proof}

\begin{rem}\label{rem:best-possible}
The hypothesis $h\le k-p^r$ of Proposition~\ref{prop:binomial} cannot be weakened in a nontrivial way.
(Of course, that hypothesis is equivalent to $h\le k+1-p^r$ when $r>0$, because of the other hypothesis $p\nmid h$.)
In fact, when $k=3p^r-1$ we have
$(-1)^i\binom{k}{i}\equiv 1\pmod{p}$
for $0\le i<p^r$ and for $2p^r\le i<3p^r$.
Therefore,
$(-1)^i\binom{k}{i}\pmod{p}$ is a periodic function of $i$ with period $h$ in the range $0\le i\le k$,
for any $h\ge 2p^r=k-p^r+1$.
For $p=3$ and $k=5p^r-1$ one checks similarly that
$(-1)^i\binom{k}{i}\pmod{p}$ is a periodic function of $i$ with period $h$ in the range $0\le i\le k$,
for any $h\ge 4p^r=k-p^r+1$.

The cases  $h=2p^r+1$ for $p\not=3$ and $h=4p^r+1$ for $p=3$
show that the hypotheses $3(h-1)<2(k+1)$ for $p\not=3$ and $5(h-1)<4(k+1)$ for $p=3$
in Theorem~\ref{thm:binomial} cannot be weakened, either.
\end{rem}

The hypothesis that $p\nmid h$ in Proposition~\ref{prop:binomial} or Theorem~\ref{thm:binomial} is also indispensable,
because if $p^t-p^s\le k<p^t$ then
$(-1)^i\binom{k}{i}\pmod{p}$ is a periodic function of $i$ with period $p^s$ in the range $0\le i\le k$.
However, Theorem~\ref{thm:binomial} can be used to prove that these values of $k$ are essentially the only exceptions
if we allow $p$ to divide $h$.

\begin{cor}\label{cor:binomial}
Let $p$ be a prime, let $h,k$ be positive integers, and let $p^s$ be the highest power of $p$ which divides $h$.
Suppose that $k\ge 5p^s$, and that either $3(h-p^s)<2(k+1)$ or $p=3$ and $5(h-p^s)<4(k+1)$.
Suppose that $(-1)^i\binom{k}{i}\pmod{p}$
is a periodic function of $i$ with period $h$ in the range $0\le i\le k$.
Then
$p^t-p^s\le k<p^t$ for some integer $t$.
\end{cor}

\begin{proof}
Express $h=h'p^s$ where $p\nmid h'$ and $k=k'p^s+k''$ where $0\le k''<p^s$.
For $i=i'p^s$ Lucas' theorem asserts that
$\binom{k}{i}\equiv\binom{k'}{i'}$.
Consequently
$(-1)^{i'}\binom{k'}{i'}$ is a periodic function of $i'$ with period $h'$ in the range $0\le i'\le k'$.
Since all hypotheses of Theorem~\ref{thm:binomial} are satisfied with $k'$ and $h'$ in place of
$k$ and $h$, we conclude that $k'+1$ is a power of $p$, and hence
$p^t-p^s\le k<p^t$ for some integer $t$.
\end{proof}

\begin{rem}
The hypotheses
$3(h-p^s)<2(k+1)$ for $p\neq 3$ and $5(h-p^s)<4(k+1)$ for $p=3$
in Corollary~\ref{cor:binomial} cannot be weakened.
This follows from the discussion in Remark~\ref{rem:best-possible}
by taking $k=3p^{r+s}-1$ and $h=2p^{r+s}+p^s$ for $p\neq 3$,
or $k=5p^{r+s}-1$ and $h=4p^{r+s}+p^s$ for $p=3$.

The hypothesis $k\ge 5p^s$ is needed for $p\neq 2,5$ because periodicity occurs with
$k=5p^s-1$ and $h=4p^s$, which are compatible with the remaining assumptions.
\end{rem}

\section{An application to multiplicative subgroups of finite fields}\label{sec:finite-fields}

The special case of Theorem~\ref{thm:binomial} where $h$ divides $k$ can be interpreted in terms of finite fields,
and yields Theorem~\ref{thm:near-field}, which we have stated in the Introduction.

\begin{proof}[Proof of Theorem~\ref{thm:near-field}]
Let $k$ be the order of $G$ and $h$ the order of $H$.
The roots of the polynomial $(1-x^k)/(1-x^h)$ are distinct and are the elements of $G\setminus H$.
By hypothesis all the roots of this polynomial are also roots of the polynomial $(1-x)^k-1$.
Hence there exists a polynomial $g(x)\in\F_p[x]$, necessarily of degree $h$
and without constant term, such that
\[
\sum_{i=1}^{k}(-1)^i\binom{k}{i}x^i=
(1-x)^k-1=g(x)\cdot(1-x^k)/(1-x^h)=g(x)\cdot\sum_{j=0}^{k/h-1}x^{jh}
\]
in $\F_p[x]$.
It follows that $(-1)^i\binom{k}{i}\pmod{p}$
is a periodic function of $i$ with period $h$ in the range $0<i\le k$.
Because of the identity $\binom{k}{k-i}=\binom{k}{i}$ the periodicity extends to the range $0\le i\le k$.
Corollary~\ref{cor:binomial-weaker} applies and yields that $k+1$ is a power of $p$.
The binomial theorem implies that $G\cup\{0\}$, the set of roots of $x^{k+1}-x$, is additively closed.
Hence $G\cup\{0\}$ is a subfield of $\F_q$
and, therefore, coincides with it.
\end{proof}

It would be interesting to know whether the general case of Theorem~\ref{thm:binomial}
(where $h$ does not necessarily divide $k$) has any useful interpretation in terms of finite fields.

An alternative proof of Theorem~\ref{thm:near-field}, which actually establishes the more general
Theorem~\ref{thm:near-field-general} below,
is based on an argument of Leep and Shapiro, which can be extracted from the proof of~\cite[Lemma~3]{LeeSha}
and stated as follows.

\begin{lemma}\label{lemma:Leep-Shapiro}
Let $F$ be any field, let $G$ be a multiplicative subgroup of $F$
and let $H$ be a finite proper subgroup of $G$.
Suppose that
$1-\alpha\in G$ for all $\alpha\in G\setminus H$.
Then $1-\alpha\in G$ for all $\alpha\in G\setminus \{1\}$.
\end{lemma}

\begin{proof}
Arguing by contradiction we assume that there exists $\beta\in H\setminus\{1\}$ with $1-\beta\not\in G$.
Then for every $\alpha\in G\setminus H$ we have
$1-\alpha\in G$ and $1-\alpha\beta\in G$, because $\alpha\beta\in G\setminus H$.
Since $(1-\alpha\beta)-(1-\alpha)=\alpha(1-\beta)\not\in G$,
we know that
$\tilde\alpha=(1-\alpha)/(1-\alpha\beta)\in G$
satisfies $1-\tilde\alpha\not\in G$,
and hence $\tilde\alpha\in H$.
Also, $\tilde\alpha\not=1/\beta$ because $\beta\not=1$.
Thus, the map $\alpha\mapsto\tilde\alpha$ sends $G\setminus H$, injectively, into $H\setminus\{1/\beta\}$.
Since $H$ is a proper finite subgroup of $G$ and $G\setminus H$ is a union of cosets of $H$ we reach a contradiction.
\end{proof}

\begin{theorem}\label{thm:near-field-general}
Let $F$ be any field, let $G$ be a multiplicative subgroup of $F$
and let $H$ be a finite proper subgroup of $G$.
If
$1-\alpha\in G$ for all $\alpha\in G\setminus H$,
then $G\cup\{0\}$ is a subfield.
\end{theorem}

Lemma~\ref{lemma:Leep-Shapiro} reduces the proof of Theorem~\ref{thm:near-field-general}
to the case $H=1$, and then it follows that $G\cup\{0\}$ is a subfield of $F$
by an argument which we have given after stating Theorem~\ref{thm:binomial}.
Note that both Lemma~\ref{lemma:Leep-Shapiro} and Theorem~\ref{thm:near-field-general}
would be false without the finiteness assumption on $H$.
A counterexample (taken from the proof of~\cite[Proposition~6]{LeeSha})
is obtained by taking $F=\Q_p$, $G=v_p^{-1}(m\Z)$ for some integer $m>1$,
and $H=1+p\Z_p$.

\medskip

Writing $1-G=\{1-g:g\in G\}$ we can restate Theorem~\ref{thm:near-field} as follows:
if $G$ is a subgroup of $\F_q^\ast$ such that
$G\setminus H\subseteq G\cap(1-G)$
for some proper subgroup $H$ of $G$, then $G\cup\{0\}$ is a field.
The set $G\cap(1-G)$ which appears in the hypothesis is closely related with the solutions
of the equation $x^n+y^n=1$, where $n=(q-1)/k$.
In fact, since $G=\{\beta^n:\beta\in\F_q^\ast\}$,
the set of solutions of $x^n+y^n=1$ in $\F_q^\ast\times\F_q^\ast$
is in a $n^2$-to-one correspondence with $G\cap(1-G)$, given by $(\beta,\gamma)\mapsto\beta^n$.
In particular,
$|G\cap(1-G)|=(N-d)/n^2$,
where $N$ is the number of projective $\F_q$-rational points of the Fermat curve $x^n+y^n=z^n$
(written in homogeneous coordinates), and $d$ is the number of such points with $xyz=0$.
Clearly, $d=3n$ if $|G|$ is even and $d=2n$ otherwise.
Since the hypothesis that
$G\setminus H\subseteq G\cap(1-G)$
in the formulation of Theorem~\ref{thm:near-field} given above implies that
$|G\cap(1-G)|\ge |G|/2$,
it is natural to ask whether the conclusion follows from this weaker condition
via known bounds for $N$.
The following result shows how far one can get using
Weil's bound $|N-q-1|\le (n-1)(n-2)\sqrt{q}$
(see~\cite{IR} or~\cite{LN}).

\begin{theorem}\label{thm:near-field-Weil}
Let $G$ be a multiplicative subgroup of the finite field $\F_{q}$,
and suppose that $|G\cap(1-G)|\ge |G|/2$.
Then either $G=\F_q^\ast$ or $|G|<2\sqrt{q}$.
\end{theorem}

\begin{proof}
Let $k=|G|$ and $n=(q-1)/k$ as above.
Recalling that $d\ge 2n\ge 2$, Weil's upper bound implies the weaker inequality
$N\le q-1+n^2\sqrt{q}+d$, and hence
$|G\cap(1-G)|\le k^2/(q-1)+\sqrt{q}$.
Because of our hypothesis it follows that $f(k)\ge 0$, where
$f(X)=X^2-\frac{q-1}{2}X-(q-1)\sqrt{q}$.
We may assume $q>5$, because the remaining cases can be easily dealt with individually.
Then we have $f((q-1)/2)=-(q-1)\sqrt{q}<0$ and
$f(2\sqrt{q})=-2\sqrt{q}(q-1-2\sqrt{q})<0$.
We conclude that either $k>(q-1)/2$ or $k<2\sqrt{q}$.
Because $k$ is a proper divisor of $q-1$, the former case yields that $k=q-1$.
\end{proof}

\begin{rem}
The following more precise form of Theorem~\ref{thm:near-field-Weil} can be
proved by a more careful application of Weil's bound, which we omit for brevity:
if $G$ is a subgroup of $\F_q^\ast$ with $|G\cap(1-G)|=c|G|$ for some $c\ge 1/2$, then either $G=\F_q^\ast$
or $|G|<\frac{r(r-1)}{c(r+1)}$, where $r=\sqrt{q}$.
This stronger formulation is seen to be best possible by taking $q$ to be a square and $G=\F_r^\ast$,
in which case
$c=(r-2)/(r-1)$, and hence
$\frac{r(r-1)}{a(r+1)}-|G|=2\frac{r-1}{(r-2)(r+1)}$
can be made arbitrarily small.
\end{rem}

A limitation of Weil's bound is that it is far from optimal, and eventually becomes trivial,
when $n$ is large
with respect to $q$.
Garc\'{\i}a and Voloch proved in~\cite{GarVol} a series of bounds
for the number $N$ of $\F_q$-rational projective points of the Fermat curve $ax^n+by^n=z^n$,
where $a,b$ are nonzero elements of $\F_q$,
which are better than Weil's bound for $n$ relatively large with respect to $q$.
Assume that $n$ is prime to $p$ and let $d$ denote the number of $\F_q$-rational points of the curve with $xyz=0$.
Then the first of Garc\'{\i}a and Voloch's bounds reads
\begin{equation*}
N\le n(n+q-1-d)/2+d.
\end{equation*}
According to~\cite[Theorem~2]{GarVol}, this bound
holds for $q$ odd except when $a,b\in\F_{p^t}$ for some subfield $\F_{p^t}$ of $\F_q$ and $n=(q-1)/(p^t-1)$.
This bound is better than Weil's bound, roughly, when $n\ge \sqrt{q}/2$.
Note that this range is roughly what we need in order to exclude the second alternative conclusion of
Theorem~\ref{thm:near-field-Weil}, under the additional assumption that
$G$ is not contained in a proper subfield of $\F_q$.
Proceeding more formally,
in the special case where $a=b=1$ and $n$ divides $q-1$, Garc\'{\i}a and Voloch's bound is equivalent
to the following statement.

\begin{theorem}\label{thm:near-field-GV}
Let $G$ be a multiplicative subgroup of order $k$ of the finite field $\F_{q}$ of odd characteristic.
Suppose $G\cup\{0\}$ is not a field.
Then $|G\cap(1-G)|\le (|G|-1)/2$.
\end{theorem}

\begin{proof}
Setting $n=(q-1)/k$, the hypothesis that $G\cup\{0\}$ is not a field
insures that the condition on $n$ for the validity of Garc\'{\i}a and Voloch's bound is met.
We have noted earlier that $d=3n$ if $k$ is even and $d=2n$ if $k$ is odd.
In terms of
$|G\cap(1-G)|=(N-d)/n^2$,
Garc\'{\i}a and Voloch's bound becomes
$|G\cap(1-G)|\le \frac{k+1}{2}-\frac{d}{2n}$,
and the latter quantity equals the integral part of $(|G|-1)/2$.
\end{proof}

Thus, the case of Theorem~\ref{thm:near-field} where $q$ is odd is also a consequence of Theorem~\ref{thm:near-field-GV}.

\section{Variations}\label{sec:variations}

We offer an analogue of Theorem~\ref{thm:binomial} for {\em unsigned} binomial
coefficients $\binom{k}{i}$.
We may assume that both $p$ and $h$ are odd, otherwise
Theorem~\ref{thm:binomial} applies.
If we keep the rest of the hypotheses of Theorem~\ref{thm:binomial} unchanged
for convenience, these additional assumptions
exclude at once that $k+1$ is a power of $p$, because $\binom{p^r-1}{i}\pmod{p}$
is periodic in the range $0\le i\le p^r-1$ with minimum period two.
However, $\binom{2p^r-1}{i}\pmod{p}$ is periodic in the range $0\le i\le 2p^r-1$ with period
any odd $h$ with $p^r\le h\le 2p^r-1$, and
$\binom{4\cdot 3^r-1}{i}\pmod{3}$ is periodic in the range $0\le i\le 4\cdot 3^r-1$ with period
any odd $h$ with $3\cdot 3^r\le h\le 4\cdot 3^r-1$.
(No other odd $h$ is a period in either case.)
Thus, an analogue of Theorem~\ref{thm:binomial} for unsigned binomial coefficients
which can be obtained without extra effort is as follows.
We need only note that the analogue of Proposition~\ref{prop:binomial} holds and is proved in the same way,
after the signed binomial coefficients $(-1)^i\binom{k}{i}$
are replaced with unsigned binomial coefficients $\binom{k}{i}$.

\begin{theorem}\label{thm:binomial-unsigned}
Let $p$ be an odd prime, and let $h,k$ be positive integers with $h$ odd and $p\nmid h$.
Suppose that $k\ge 5$, and that either $3h<2k+5$ or $p=3$ and $5h<4k+9$.
Suppose that the binomial coefficient $\binom{k}{i}$ viewed modulo $p$
is a periodic function of $i$ with period $h$ in the range $0\le i\le k$.
Then $k+1$ is twice or four times a power of $p$, with $p=3$ in the latter case,
or $(p,k,h)=(3,7,7)$.
\end{theorem}

\begin{cor}\label{cor:binomial-unsigned-weaker}
Let $p$ be a prime, and let $h,k$ be positive integers with $p\nmid h$ and $2h\le k$.
Suppose that the binomial coefficient $\binom{k}{i}$ viewed modulo $p$
is a periodic function of $i$ with period $h$ in the range $0\le i\le k$.
Then $k+1$ is a power of $p$.
\end{cor}

\begin{proof}
We may assume $p$ and $h$ odd, otherwise Corollary~\ref{cor:binomial-weaker} applies.
After inspection of the cases $k=2,3,4$ we may assume $k\ge 5$.
Theorem~\ref{thm:binomial-unsigned} yields that either
$k=2p^r-1$, or $p=3$ and $k=4\cdot 3^r-1$, for some $r$.
The discussion preceding the theorem shows that
the minimum period $h$ prime to $p$ equals $p^r+1$ or $3\cdot 3^3+1$, respectively.
In both cases $2h>k$, contradicting one of our hypotheses.
\end{proof}

The following application of Corollary~\ref{cor:binomial-unsigned-weaker} follows
in the same way as Theorem~\ref{thm:near-field} follows from Corollary~\ref{cor:binomial-weaker}.

\begin{theorem}\label{thm:near-field-unsigned}
Let $G$ be a multiplicative subgroup of a finite field $\F_q$
which generates $\F_q$ as a field,
and let $H$ be a proper subgroup of $G$.
Suppose that
$\alpha+1\in G$ for all $\alpha\in G\setminus H$.
Then $G=\F_{q}^\ast$.
\end{theorem}

As we mentioned in Section~\ref{sec:intro},
when dealing with binomial coefficients modulo a prime it is
sometimes convenient to endow them with suitable alternating signs,
and to consider, for example, the integers
$(-1)^{i}\binom{k}{i}$ (the {\em signed binomial coefficients} used so far in this paper)
or $(-1)^{k}\binom{k}{i}$.
One reason in favour of the latter is that the basic recursion
$\binom{k+1}{i+1}=\binom{k}{i}+\binom{k}{i+1}$
satisfied by the binomial coefficients then takes the more symmetric form
\[
(-1)^{k}\binom{k}{i}+(-1)^{k}\binom{k}{i+1}+(-1)^{k+1}\binom{k+1}{i+1}=0.
\]
An immediate consequence of this fact
is that a version of ``Pascal's triangle'' modulo $p$
endowed with alternating signs as above and suitably truncated at $k<p^s$
displays a symmetry group isomorphic with $S_3$ (the symmetric group on three objects),
as opposed to the only symmetry given by $\binom{k}{i}=\binom{k}{k-i}$
in characteristic zero.
This is because such a symmetry is displayed by both the above recursion rule
and the three ``boundary conditions''
$\binom{k}{-1}=\binom{k}{k+1}=0$ for $k\ge 0$ and
$\binom{p^s}{i}\equiv 0\pmod{p}$ for $0<i<p^s$,
which together determine the binomial coefficients uniquely.
More formally stated, the symmetry group under consideration is generated, together with the
ordinary identity
$\binom{k}{i}=\binom{k}{k-i}$ for $k\ge 0$,
by the modular identity
\[
(-1)^{k}\binom{k}{i}\equiv
(-1)^{i}\binom{p^s-1-i}{p^s-1-k}\pmod{p}
\qquad\textrm{for $0\le i\le k\le p^s-1$.}
\]

Alternatively, the latter can be proved by combining the characteristic zero identities
$\binom{k}{i}=\binom{k}{k-i}$ and
\[
\binom{-k}{i}=(-1)^{i}\binom{k+i-1}{i}
\]
(which is valid for arbitrary integers $k$ and $i$)
with the periodicity in the upper entry with period $p^s$ mentioned in Section~\ref{sec:intro}, as follows:
\begin{align*}
(-1)^{k}\binom{k}{i}=
(-1)^{k}\binom{k}{k-i}=
(-1)^{i}\binom{-i-1}{k-i}
&\equiv
(-1)^{i}\binom{p^s-i-1}{k-i}\pmod{p}
\\
&=(-1)^{i}\binom{p^s-1-i}{p^s-1-k}.
\end{align*}

This identity allows one to translate back and forth between modular properties
of (portions of) rows and columns of Pascal's triangle,
provided suitable alternating signs are introduced.
In particular, we deduce at once from Theorem~\ref{thm:binomial}
an analogous result concerning periodicity of $\binom{k}{i}\pmod{p}$ with respect to the upper entry
in a certain range.

\begin{theorem}\label{thm:binomial-vertical}
Let $p$ be a prime, and let $h,i,s$ be positive integers with $p\nmid h$.
Suppose that $i<p^s-5$, and that either $3h<2p^s-2i+3$ or $p=3$
and $5h<4p^s-4i+5$.
Suppose that the binomial coefficient $\binom{k}{i}$ viewed modulo $p$
is a periodic function of $k$ with period $h$ in the range $i\le k\le p^s-1$.
Then $i+1$ is a power of $p$.
\end{theorem}

\bibliography{References}

\end{document}